\documentclass[12pt]{article}
\usepackage[russian]{babel}
\usepackage{latexsym}
\usepackage{amsmath}
\usepackage{amssymb}
\usepackage{stackrel}
\usepackage{color}
\usepackage{graphicx}
\usepackage[utf8]{inputenc}

\usepackage[left=2cm,right=2cm,top=2cm,bottom=2cm,bindingoffset=0cm]{geometry}

\begin{document}
\begin{center}\large\bf
Maillet--Malgrange type theorem for a formal Dulac series solution of an analytic ODE
\end{center}

\begin{center}\bf
Goryuchkina Irina
\end{center}

\begin{quote}\small
A Maillet-Malgrange type theorem is proved for a Dulac series (in the general case, with complex exponents), which formally satisfies an analytical ordinary differential equation (ODE). This theorem allows to estimate the growth of the norms of the coefficients of such a series, that is, to determine its Gevrey order, and in the special case it provides a sufficient condition for the convergence of the series.
\end{quote}
\bigskip\bigskip

\begin{center}
\large\bf Теорема типа Майе--Мальгранжа для ряда Дюлака, формально удовлетворяющего аналитическому ОДУ
\end{center}
\begin{center}
 \bf Горючкина И.В.
\end{center}

\begin{quote}\small
Доказывается теорема типа Майе--Мальгранжа для ряда Дюлака (в общем случае, с комплексными показателями степени), который формально удовлетворяет {\it аналитическому} обыкновенному дифференциальному уравнению (ОДУ). Данная теорема позволяет оценить скорость роста норм коэффициентов такого ряда, то есть определить его порядок Жевре, а в частном случае представляет достаточное условие сходимости ряда.
\end{quote}

\bigskip
Ряд Дюлака --- это формальный ряд вида
\begin{equation}\label{eq2}
  \sum\limits_{k=1}^{\infty}c_k(t)\,x^{\lambda_k},\qquad t:=\ln x, \qquad c_k\in\mathbb{C}[t],\qquad c_1(t)\not\equiv 0,
\end{equation}
$$
\lambda_k\in\mathbb{C}, \qquad 0<{\rm Re}\,\lambda_1\leqslant{\rm Re}\,\lambda_2\leqslant\ldots,\qquad {\rm Re}\,\lambda_k\rightarrow\infty.
$$
Нас будут интересовать ряды вида \eqref{eq2}, которые являются формальными решениями {\it аналитического} ОДУ\footnote{Позже мы обсудим корректность подстановки формального ряда \eqref{eq2} в уравнение \eqref{eq1} и определим понятие формального решения.} произвольного порядка $n$,
\begin{equation}\label{eq1}
F(x,y,\delta y,\ldots,\delta^ny)=0,\qquad \delta=x\frac{d}{dx},
\end{equation}
где   $F(x,y_0,y_1,\ldots,y_n)\not\equiv 0$ --- голоморфная вблизи $0\in{\mathbb C}^{n+2}$ функция, $F(0,\ldots,0)=0$.

\bigskip
Впервые ряды вида \eqref{eq2} при $\lambda_k\in\mathbb R$ появились в работе Анри Дюлака \cite{Du} по предельным циклам векторного поля на плоскости
как асимптотические разложения отображения монодромии в окрестности гиперболического полицикла.
В настоящий момент ряды Дюлака встречаются в разных областях математики, в том числе, в теории динамических систем, а также среди формальных решений 
 аналитических ОДУ. Первые шаги в изучении <<неформальных>> свойств таких формальных решений были сделаны в \cite{GG_D}, где были предложены достаточные условия сходимости формального ряда Дюлака с $\lambda_k=k\in\mathbb{Z}_+$, который удовлетворяет {\it алгебраическому\,} ОДУ. Здесь мы обобщаем  результаты работы \cite{GG_D}. А именно, доказываем, что если формальный ряд \eqref{eq2} удовлетворяет {\it аналитическому\,} ОДУ \eqref{eq1}, то существуют числа $C,$ $A>0$ и значение $s\in\mathbb{R}_{>0}\cup \{+\infty\}$ (которое определяется по уравнению \eqref{eq1} и ряду \eqref{eq2}, см. теорему 1 далее) такие, что
\begin{equation}\label{dd1}
\|c_k\|_R\leqslant C A^k\, \left|\Gamma\left(\lambda_k/s\right)\right|,
\end{equation}
где $\Gamma$ --- это гамма-функция Эйлера и норма мероморфной в окрестности бесконечности функции $f(t)=\sum\limits_{i\geqslant p}f_i\; t^{-i}$ определяется как
$$
\|f\|_R=\sum\limits_{i\geqslant p}|f_i|\, R^{-i},
$$
при этом $R>1$ --- фиксированная величина, определяемая далее.
\medskip

Число $1/s$ назовём {\it порядком Жевре} ряда \eqref{eq2}, удовлетворяющего условию (\ref{dd1}).
\medskip

Отметим, что в случае $s=+\infty$ имеется равномерная сходимость ряда \eqref{eq2} в каждом секторе с вершиной в нуле достаточно малого радиуса и раствора меньше $2\pi$.\footnote{В таком случае сходимость вызвана не только скоростью роста нормы коэффициентов ряда, но и свойствами множества показателей степени ряда, формально удовлетворяющего ОДУ, --- см. теорему 2 далее.}

В случае  $s<+\infty$  оценка \eqref{dd1} является первым шагом на пути исследования вопроса существования в каждом открытом секторе $V$ с вершиной в нуле достаточно малых радиуса и раствора решения $f_V(x)$ уравнения \eqref{eq1} такого, что для любого подсектора $W\subset V$ существуют положительные числа $C=C(W)$ и $A=A(W)$ такие, что выполняется неравенство $\displaystyle\Bigl|f_V(x)-\sum\limits_{k=1}^{N-1}c_k(t)\,x^{\lambda_k} \Bigr|\leqslant C A^{N}\, \left|\Gamma\left(\lambda_N/s\right)\right| \left|x^{\lambda_N}\right|$ при всех $x\in W$, $N\in\mathbb{N}$. Для формальных рядов Тейлора удовлетворяющих уравнению \eqref{eq1} это было сделано в работе \cite{RamisSibuya}.

\bigskip
Основной результат этой работы можно сформулировать следующим образом.

\bigskip
{\bf Теорема 1}. {\it Пусть ряд \eqref{eq2} формально удовлетворяет уравнению \eqref{eq1}, выполнено условие $\displaystyle F'_{y_n}(x,\varphi,\dots,\delta^n \varphi)\not\equiv 0$ и для каждого $\,j=0,\,1,\,\ldots,\,n\,$ имеем ряд Дюлака
\begin{eqnarray*}
F'_{y_j}(x,\varphi,\dots,\delta^n \varphi)=A_j\,x^{\nu}+B_j(t)\, x^{\nu_j}+\ldots, &\quad& A_j\in{\mathbb C}, \quad
B_j\in{\mathbb C}[t]\setminus\{0\},  \\ &  & {\rm Re}\,\nu_j>{\rm Re}\,\nu,
\end{eqnarray*}
где $\nu\in\mathbb{C}$ --- одинаковое для всех $\,j=0,\,1,\,\ldots,\,n$,\, и хотя бы один из коэффициентов $A_j$ отличен от нуля.
Положим $\;\ell=\max\{j\;|\;A_j\neq 0\}\;$ и
\begin{eqnarray}\label{slope}
s=\left\{\begin{array}{rl}
            +\infty, & \mbox{ если } A_n \neq 0,\\
            \\
           \displaystyle\min_{j>\ell}\frac{{\rm Re}\,\nu_j-{\rm Re}\,\nu}{j-\ell}, & \mbox{ если } A_n=0.
         \end{array}\right.
\end{eqnarray}
Тогда для всякого открытого сектора $S$ с вершиной в нуле достаточно малого радиуса и раствора меньше $2\pi$ ряд
$$
\sum_{k=1}^{\infty}\frac{c_k(t)}{\Gamma(\lambda_k/s)}\,x^{\lambda_k}
$$
сходится равномерно в $S$.}

\medskip
Теорема~1 обобщает известную теорему Майе--Мальгранжа \cite{Mal} для формального {\it степенного} ряда, удовлетворяющего аналитическому ОДУ (\ref{eq1}). В наших предыдущих работах \cite{GG_MM} и \cite{GG_D} рассматривались частные случаи теоремы 1 для ряда \eqref{eq2} и {\it алгебраического} ОДУ (\ref{eq1}): в \cite{GG_MM} коэффициенты $c_k\in\mathbb{C}$ постоянны, а в \cite{GG_D} показатели степени $\lambda_k=k\in\mathbb{Z}_+$ и $s=+\infty$.

\medskip
Далее обсудим некоторые свойства рядов \eqref{eq2}, удовлетворяющих уравнению \eqref{eq1}. А именно: корректность подстановки ряда \eqref{eq2} в уравнение \eqref{eq1}, принадлежность показателей степени $\lambda_k$ к конечно порожденной аддитивной полугруппе с линейно независимыми над $\,\mathbb Z\,$ образующими и возможность перехода от ряда \eqref{eq2} к ряду Тейлора многих переменных.

\bigskip
Заметим, что множество всех формальных рядов вида \eqref{eq2}
$$
\mathcal{D}=\Bigl\{\varphi=\sum\limits_{k=1}^{+\infty} c_k(t)\, x^{\lambda_k} \mid \quad c_k\in\mathbb{C}[t],\quad  \lambda_k\in{\mathbb C}, \quad 0<{\rm Re}\,\lambda_1\leqslant{\rm Re}\,\lambda_2\leqslant\ldots\rightarrow\infty\Bigr\}
$$
является кольцом. Действительно, пусть
$$
\varphi=\sum\limits_{k=1}^{+\infty} c_k(t)\,x^{\lambda_k}\in\mathcal{D}, \qquad \tilde\varphi=\sum\limits_{k=1}^{+\infty}\tilde c_k(t)\,x^{\tilde\lambda_k}\in\mathcal{D},
$$
тогда $\varphi+\tilde{\varphi}\in\mathcal{D}$, так как объединение множеств показателей степени таких рядов снова может быть упорядочено по неубыванию вещественных частей показателей, а произведение двух рядов
представляется в виде
$$
\varphi\,\tilde\varphi=\sum_{\sigma\in\mathbb C}\Bigl(\sum\limits_{\lambda_i+\tilde\lambda_j=\sigma} c_i(t)\,\tilde c_j(t)\Bigr)x^\sigma,
$$
при этом внутренние суммы корректно определены, поскольку для каждого $\sigma$ найдётся только конечное число пар $(\lambda_i,\tilde\lambda_j)$ таких, что $\lambda_i+\tilde\lambda_j=\sigma$, в силу стремления последовательностей $\{{\rm Re}\,\lambda_k\}_{k\in\mathbb N}$, $\{{\rm Re}\,\tilde\lambda_k\}_{k\in\mathbb N}$ к бесконечности. По этой же причине для каждого $N>0$ имеется только конечное число пар $(\lambda_i,\tilde\lambda_j)$ таких, что ${\rm Re}\,\lambda_i+{\rm Re}\,\tilde\lambda_j<N$, то есть в ряду $\varphi\,\tilde\varphi$ имеется только конечное число показателей $\sigma$ таких, что ${\rm Re}\,\sigma<N$ и, следовательно, данный ряд также является рядом Дюлака.

\medskip
Действующий на пространстве голоморфных в окрестности нуля функций дифференциальный оператор $\delta=x(d/dx): y(x)\mapsto xy'(x)$ естественным образом  продолжается на кольцо формальных рядов Дюлака $\cal D$:
$$
\delta\left(\sum_{k=1}^{+\infty}c_k(t)\,x^{\lambda_k}\right)=\sum_{k=1}^{+\infty}\left(\lambda_k+\frac{d}{dt}\right)c_k(t)\, x^{\lambda_k}.
$$
Таким образом, если $F(x,y_0,y_1,\ldots,y_n)$ --- многочлен и $\varphi\in\cal D$, то $F(x,\varphi,\delta\varphi,\ldots,\delta^n\varphi)\in\cal D$. Для корректности подстановки ряда \eqref{eq2} не только в полиномиальное, но и в {\it аналитическое} ОДУ, нам нужно показать, что то же самое справедливо и в случае если $F$ --- голоморфная функция в окрестности $0\in{\mathbb C}^{n+2}$.
\bigskip

{\bf Лемма 1.}  {\it Пусть $F(x,y_0,y_1,\ldots,y_n)$ --- голоморфная функция в окрестности $0\in{\mathbb C}^{n+2}$, $F(0)=0$ и $\varphi\in\cal D$. Тогда $F(x,\varphi,\delta\varphi,\ldots,\delta^n\varphi)\in\cal D$.}
\medskip

{Д о к а з а т е л ь с т в о } повторяет рассуждения леммы 1 из \cite{GGarxiv}. Представим функцию $F$ в окрестности $0\in{\mathbb C}^{n+2}$ в виде ряда
$$
F(x,y_0,y_1,\ldots,y_n)=F_1(x,y_0,y_1,\ldots,y_n)+F_2(x,y_0,y_1,\ldots,y_n)+\ldots,
$$
где $F_d$ --- однородный многочлен степени $d$. Тогда показатели степени каждого ряда Дюлака $F_d(x,\varphi,\delta\varphi,\ldots,\delta^n\varphi)\in\cal D$,
$$
F_d(x,\varphi,\delta\varphi,\ldots,\delta^n\varphi)=\sum_{k=1}^{+\infty}c_k^{(d)}(t)\,x^{\lambda_k^{(d)}}, \qquad d=1,2,\ldots,
$$
обладают следующим свойством: $0<d\min(1,{\rm Re}\,\lambda_1)\leqslant{\rm Re}\,\lambda_1^{(d)}\leqslant{\rm Re}\,\lambda_2^{(d)}\leqslant\ldots\rightarrow\infty$, где $\lambda_1$ --- первый показатель степени ряда $\varphi$, то есть для каждого $N>0$ лишь конечное число показателей степени $\lambda_k^{(d)}$ в сумме $\displaystyle\sum\limits_{d=1}^{+\infty}F_d(x,\varphi,\delta\varphi,\ldots,\delta^n\varphi)$ удовлетворяют условию ${\rm Re}\,\lambda_k^{(d)}<N$ (а именно, это условие может быть выполнено только при $d<N/\min(1,{\rm Re}\,\lambda_1)$, и при каждом таком $d$ имеется лишь конечное число $\lambda_k^{(d)}$, для которых ${\rm Re}\,\lambda_k^{(d)}<N$, поскольку $F_d(x,\varphi,\delta\varphi,\ldots,\delta^n\varphi)\in\cal D$) и, тем самым, $F(x,\varphi,\delta\varphi,\ldots,\delta^n\varphi)\in\cal D$. {\hfill$\Box$}
\bigskip

Доказанная лемма объясняет корректность следующего определения формального решения \eqref{eq2} уравнения \eqref{eq1}.
\bigskip

Будем говорить, что ряд Дюлака $\varphi\in\cal D$ является {\it формальным решением} уравнения \eqref{eq1}, если $F(x,\varphi,\delta\varphi,\ldots,\delta^n\varphi)$ --- это ряд Дюлака с нулевыми коэффициентами:
$$
F(x,\varphi,\delta\varphi,\ldots,\delta^n\varphi)=0\in{\cal D}.
$$

Жирным шрифтом будем обозначать целочисленные мультииндексы и другие векторные величины, при этом их размерность иногда указывается явно, а иногда понятна из контекста. Модуль мультииндекса, как обычно, обозначает сумму координат последнего. Например, $|{\bf p}|=p_1+\ldots+p_n$ для мультииндекса ${\bf p}=(p_1,\ldots, p_n)$.
\medskip

Запишем формальный ряд \eqref{eq2} в виде
\begin{equation}\label{eq3}
\varphi=\sum\limits_{k=1}^{m}c_k(t)\,x^{\lambda_k}+x^{\lambda_m}\psi=\varphi_{m}+x^{\lambda_m}\psi,
\end{equation}
где
\begin{equation}\label{eq4}
\psi=\sum_{k=m+1}^{\infty}c_k(t)\,x^{\lambda_k-\lambda_m}.
\end{equation}


Обозначим $\mathcal{D}^\circ\subset \mathcal{D}$ кольцо рядов \eqref{eq2}, сходящихся в некотором секторе $S$ с вершиной в нуле достаточно малого радиуса и раствора меньше $2\pi$. Кроме того, пусть $${\boldsymbol\varphi}=(\varphi,\delta \varphi,\ldots,\delta^n\varphi),\qquad{\boldsymbol\varphi}_m=(\varphi_m,\delta \varphi_m,\ldots,\delta^n\varphi_m), \qquad{\boldsymbol\psi}=(\psi,(\lambda_m+\delta)\psi,\ldots,(\lambda_m+\delta)^n\psi).$$
Также определим понятие {\it порядка} ${\rm val}\,\varphi$ ряда Дюлака следующим образом:
$${\rm val}\,\varphi={\rm Re}\,\lambda_1.$$

\bigskip
{\bf Лемма 2.} {\it В условиях теоремы $1$ найдется такой номер $m\in\mathbb{N}$, что уравнение \eqref{eq1} с помощью преобразования
\begin{equation}\label{eq5}
y=\varphi_{m}+x^{\lambda_m}u
\end{equation}
приводится к уравнению специального вида
\begin{eqnarray}\label{reduced}
L(\lambda_m+\delta)u+ \widetilde{L}(t,\,x,\,\lambda_m+\delta)u+N(t,\,x,\,x^\tau u,\,x^\tau(\lambda_m+\delta) u,\,\dots,\,x^\tau(\lambda_m+\delta)^n u)=0,
\end{eqnarray}
где $\,L(\lambda_m+\delta)\,$ и $\,\widetilde{L}(t,\,x,\,\lambda_m+\delta)\,$ --- это линейные дифференциальные операторы,
$$
L(\zeta)=\sum_{j=0}^{\ell}A_j\, \zeta^j,\qquad L(\lambda_m+z)\neq 0 \qquad \mbox{ при }\qquad {\rm Re}\,z\geqslant0,
$$
$$
\widetilde{L}(t,\,x,\,\zeta)=\sum_{j=0}^{n}\bigl(B_j(t)\,x^{\mu_j}+\ldots\bigr)\zeta^j\in\mathcal{D}^{\circ}[\zeta], \quad \mu_j=\nu_j-\nu,
$$
$$
N(t,\,x,\,u_0,\,u_1,\,\dots,\,u_n)=\sum\limits_{{\bf q}\in\mathbb{Z}_+^{n+1}}a_{\bf q}(t,\,x)\;u_0^{q_0}\ldots u_n^{q_n},\qquad
a_{\bf q}(t,\,x)\in\mathcal{D}^{\circ},
$$
$N(t,\,x,\,u_0,\,u_1,\,\dots,\,u_n)$ --- это голоморфная в окрестности $\;0\in\mathbb{C}^{n+1}\;$ функция переменных $u_0,\dots, u_n$, не содержащая их линейные степени, $\tau=(n-\ell)s$ при $\ell<n$ $(=0$ при  $\ell=n)$.}

\medskip
{Д о к а з а т е л ь с т в о.}
Разложим в (формальный) ряд Тейлора левую часть уравнения $F(x,\;{\boldsymbol \varphi})=0$. Получим равенство
\begin{eqnarray}\label{Taylor}
0 & = & F(x, {\boldsymbol\varphi}_m+x^{\lambda_m}{\boldsymbol\psi})=F(x,\, {\boldsymbol\varphi}_m)+x^{\lambda_m}\sum_{i=0}^n\frac{\partial F}{\partial y_i}(x,\,{\boldsymbol\varphi}_m)\psi_i+\nonumber \\
  &  & +\frac12x^{2\lambda_m}\sum_{i,j=0}^n\frac{\partial^2F}{\partial y_i\partial y_j}(x,\,{\boldsymbol\varphi}_m)\psi_i\psi_j+\ldots,
\end{eqnarray}
где $\psi_i=(\lambda_m+\delta)^i\psi$.


\bigskip
Выберем $m$ так, чтобы выполнялись неравенства
\begin{equation}\label{d3}
L(\lambda_m+z)\neq 0 \; \mbox{ при }\;{\rm Re}\,z\geqslant 0,
\end{equation}
$$
{\rm Re}(\lambda_{m+1}-\lambda_m)>0,\qquad {\rm Re}\,\lambda_m>\max\limits_{j=0,\,1,\,\ldots,\,n}{\rm Re}\,\nu_j+2\tau.
$$
Тогда для каждого $i=0,1,\ldots,n$ будем иметь
$$
{\rm val}\left(\frac{\partial F}{\partial y_i}(x,\boldsymbol \varphi)-\frac{\partial F}{\partial y_i}(x,\,\boldsymbol\varphi_m)\right)=
{\rm val}\biggl(x^{\lambda_m}\sum_{j=0}^n\frac{\partial^2 F}{\partial y_i\partial y_j}(x,\,\boldsymbol\varphi_m)\psi_j+\ldots\biggr)>{\rm Re}\,\lambda_m
>{\rm Re}\,\nu_i+2\tau.
$$
Из этого следует, что
\begin{eqnarray}\label{decomposition}
\frac{\partial F}{\partial y_i}(x,\, {\boldsymbol\varphi}_m)=A_i\,x^{\nu}+B_i(t)\,x^{\nu_i}+\ldots,
\end{eqnarray}
откуда
$$
x^{-\nu}\sum_{i=0}^n\frac{\partial F}{\partial y_i}(x,\,\boldsymbol\varphi_m)(\lambda_m+\delta)^i=L(\lambda_m+\delta)+\widetilde{L}(t,\,x,\,\lambda_m+\delta).
$$
Из  \eqref{Taylor} и \eqref{decomposition} следует
$$
{\rm val}\,F(x,\, {\boldsymbol\varphi}_m)>{\rm Re}(\lambda_m+\nu).
$$
Тогда, разделив соотношение \eqref{Taylor} на $x^{\lambda_m+\nu}$, мы получим, что  ряд $\psi$ \eqref{eq4} удовлетворяет уравнению \eqref{reduced}.   $\hfill\Box$
\bigskip


{\bf Замечание 1.} Из условия теоремы 1 и определения (\ref{slope}) величины $s$ следует, что
$$
{\rm Re}\,\mu_j>0, \qquad {\rm Re}\,\mu_j\geqslant(j-\ell)s, \quad j=0,\,1,\,\ldots,\,n.
$$

{\bf Замечание 2.} Коэффициенты $c_k(t)$ формального решения \eqref{eq4} полученного в лемме 2 уравнения \eqref{reduced}  определяются полиномами $c_1,$ $\ldots,$ $c_m$ однозначно. Действительно, каждый коэффициент $c_k(t)$ с $k\geqslant m+1$ удовлетворяет неоднородному линейному ОДУ с постоянными коэффициентами:
\begin{equation}\label{coef_eq}
L\Bigl(\lambda_k+\frac{d}{dt}\Bigr)v(t)=b_k(t),
\end{equation}
где многочлен $b_k(t)$ однозначно определяется предыдущими коэффициентами $c_{k-1}(t),\ldots,c_1(t)$. Поскольку выполняется условие \eqref{d3}, то есть $z=0$ не является корнем характеристического уравнения $L\left(\lambda_k+z\right)=0$, то уравнение \eqref{coef_eq} имеет единственное решение в виде многочлена.

\medskip
Отметим, что пользуясь методами из \cite{b}, можно вычислить сумму $\varphi_m$, а затем с помощью преобразования \eqref{eq5} перейти от уравнения \eqref{eq1} к уравнению \eqref{reduced}. Из замечания 2 следует, что уравнение \eqref{reduced} имеет (и при том единственное) формальное решение в виде ряда \eqref{eq4}.

Доказательство следующей вспомогательной леммы не приводим, оно аналогично доказательству леммы 4 из \cite{GGarxiv} с последующим применением леммы 5 из \cite{GGarxiv}.

\medskip
{\bf Лемма 3.} {\it В лемме $2$  ряды
$$
B_j(t)\,x^{\mu_j}+\ldots\in{\cal D}^{\circ}\qquad \mbox{и} \qquad a_{\bf q}(t,\,x)\in{\cal D}^{\circ}
$$
принадлежат кольцу $\,{\mathbb C}[t]\{x^{\rm G}\}\,$, где ${\rm G}$ --- конечно порождённая аддитивная полугруппа, образующие $r_1,\ldots,r_\kappa$ которой линейно независимы над $\mathbb Z$ и имеют положительные вещественные части.  При этом, $1$ $\lambda_1,$ $ \ldots,$ $\lambda_m,$ $\tau\in\rm G$.}
\medskip

Конечно порождённой аддитивной полугруппой ${\rm G}$ здесь мы называем множество
$$
{\rm G}=\biggl\{\sum_{j=1}^{\kappa}m_jr_j \mid  {\bf m}\in{\mathbb Z}_+^{\kappa}\setminus\{0\}\biggr\},
$$
то есть $0\not\in{\rm G}$.
Обозначение ${\mathbb C}[t]\{x^{{\rm G}}\}$ используется для кольца {\it сходящихся} (имеющих ненулевой радиус сходимости)  рядов Дюлака с показателями степени из ${\rm G}$. Обозначение ${\mathbb C}[t][[x^{{\rm G}}]]$  используется для кольца {\it формальных}  рядов Дюлака с показателями степени из ${\rm G}$.
\medskip

Лемма 3 позволяет доказать следующее (подобное теореме 1 из \cite{GGarxiv}) утверждение.

\bigskip
{\bf Теорема 2.}  {\it Пусть  ряд  $\varphi$ \eqref{eq2} удовлетворяет уравнению \eqref{eq1} и условию  теоремы $1$. Тогда
$\varphi\in{\mathbb C}[t][[x^{\rm G}]]\subset\cal D$.}
\bigskip

{Д о к а з а т е л ь с т в о.} Представим $\varphi$ в виде \eqref{eq3}, где $\psi$ удовлетворяет уравнению \eqref{reduced}, то есть
$$
L(\lambda_m+\delta)\,\psi=-\sum_{j=0}^{n}\left(B_j(t)\,x^{\mu_j}+\ldots\right) (\lambda_m+\delta)^j\psi-
$$
$$
-\sum_{{\bf q}\in{\mathbb Z}_+^{n+1}}a_{\bf q}(t,\,x){x^{\tau|{\bf q}|}}\psi^{q_0}((\lambda_m+\delta)\psi)^{q_1}\ldots ((\lambda_m+\delta)^n\psi)^{q_n},
$$
где согласно лемме 3 коэффициенты $B_j(t)\,x^{\mu_j}+\ldots,\;$  $a_{\bf q}(t,\,x)\in{\mathbb C}[t]\{x^{\rm G}\}.\;$ Следовательно, имеем соотношение
\begin{equation}\label{auxrel}
\sum_{k=m+1}^{+\infty}L\Bigl(\lambda_k+\frac{d}{dt}\Bigr)\,c_k(t)\;x^{\lambda_k-\lambda_m}=$$ $$=-\sum_{k=m+1}^{+\infty}\sum_{j=0}^{n}(B_j(t)+\ldots)\Bigl(\lambda_k+\frac{d}{dt}\Bigr)^j c_k(t)\,x^{\mu_j+\lambda_k-\lambda_m}-$$
$$
-\sum_{{\bf q}\in{\mathbb Z}_+^{n+1}}a_{\bf q}(t,\,x){x^{\tau|{\bf q}|}}\psi^{q_0}((\lambda_m+\delta)\psi)^{q_1}\ldots ((\lambda_m+\delta)^n\psi)^{q_n}.
\end{equation}
Первое слагаемое в левой части \eqref{auxrel} --- моном
$$
L\Bigl(\lambda_{m+1}+\frac{d}{dt}\Bigr)\,c_{m+1}(t)\, \,x^{\lambda_{m+1}-\lambda_m},
$$
в то время как первое слагаемое в правой части (\ref{auxrel}) --- моном, с которого начинается ряд $-a_{\bf 0}(t,\,x)\in{\mathbb C}[t]\{x^{\rm G}\}$. Поскольку выполняется условие \eqref{d3}, то $\displaystyle L\left(\lambda_{m+1}+{d}/{dt}\right)\,c_{m+1}(t)\not\equiv0$ и, следовательно, показатель степени $\lambda_{m+1}-\lambda_m$ совпадает с показателем степени монома, с которого начинается ряд $a_{\bf 0}(t,\,x)\in{\mathbb C}[t]\{x^{\rm G}\}$, то есть $\lambda_{m+1}-\lambda_m\in{\rm G}$.

Далее предположим, что $\lambda_k-\lambda_m\in{\rm G}$ при всех $k=m+1,\ldots,m+s-1$, и покажем, что тогда $\lambda_{m+s}-\lambda_m\in{\rm G}$. Число $\lambda_{m+s}-\lambda_m$ является степенью ненулевого монома
$$
L\Bigl(\lambda_{m+s}+\frac{d}{dt}\Bigr)\,c_{m+s}(t)\,x^{\lambda_{m+s}-\lambda_m},
$$
--- слагаемого с номером $s$ в левой части соотношения (\ref{auxrel}). Степень слагаемого с тем же номером (та же степень) в правой части (\ref{auxrel}) представляется в виде
$$
\alpha+\sum_{k=m+1}^{m+s-1} m_k(\lambda_k-\lambda_m),
$$
где $\alpha\in{\rm G}$, $m_k\in{\mathbb Z}_+$, то есть по индуктивному предположению является элементом полугруппы ${\rm G}$. Тем самым, $\lambda_{m+s}-\lambda_m\in{\rm G}$ и, следовательно, $\psi\in{\mathbb C}[t][[x^{\rm G}]]$, а тогда и
$\varphi\in{\mathbb C}[t][[x^{\rm G}]]$. {\hfill $\Box$}

\medskip
Линейная независимость над $\,\mathbb Z\,$ образующих $r_1,\ldots,r_{\kappa}$ аддитивной полугруппы ${\rm G}$  позволяет корректно определить биективное линейное отображение
$$
\iota:{\mathbb C}[t][[x^{{\rm G}}]]\rightarrow{\mathbb C}[t][[x_1,\ldots,x_{\kappa}]]_*
$$
из кольца ${\mathbb C}[t][[x^{{\rm G}}]]$ в кольцо ${\mathbb C}[t][[x_1,\ldots,x_{\kappa}]]_*$ формальных степенных рядов  ${\kappa}$ переменных без свободного члена с коэффициентами из ${\mathbb C}[t]$,
$$
\iota: \sum_{{\bf m}\in{\mathbb Z}_+^{\kappa}\setminus\{0\}}C_{\bf m}(t)\,x^{m_1r_1+\ldots+m_{\kappa}r_{\kappa}}\mapsto
\sum_{{\bf m}\in{\mathbb Z}_+^{\kappa}\setminus\{0\}}C_{\bf m}(t)\,x_1^{m_1}\ldots x_{\kappa}^{m_{\kappa}}.
$$
Поскольку также
$$
\iota(\eta_1\eta_2)=\iota(\eta_1)\iota(\eta_2) \qquad\forall\eta_1,\eta_2\in{\mathbb C}[t][[x^{{\rm G}}]],
$$
то отображение $\iota$ является изоморфизмом колец. Более того,  для всякого ряда вида
$$
H(x,y_0,y_1,\ldots,y_n)=\sum_{{\bf q}\in{\mathbb Z}_+^{n+1}}b_{\bf q}(t,\,x)y_0^{q_0}y_1^{q_1}\ldots y_n^{q_n}, \qquad b_{\bf q}\in{\mathbb C}[t][[x^{\rm G}]],
$$
результат подстановки формальных рядов Дюлака $\varphi_0,\varphi_1,\ldots,\varphi_n\in{\mathbb C}[t][[x^{\rm G}]]$ вместо переменных $y_0,y_1,\ldots,y_n$ в $H$ также является элементом ${\mathbb C}[t][[x^{\rm G}]]$ и имеет место следующее равенство:
\begin{equation}\label{iotagen}
\iota\bigl(H(x,\varphi_0,\varphi_1,\ldots,\varphi_n)\bigr)=\sum_{{\bf q}\in{\mathbb Z}_+^{n+1}}\iota(b_{\bf q}) \iota(\varphi_0)^{q_0}\iota(\varphi_1)^{q_1}\ldots\iota(\varphi_n)^{q_n}
\end{equation}
(доказывается аналогично лемме 6 из \cite{GGarxiv}).

Дифференцирование $\delta=x(d/dx)$, действующее на ${\mathbb C}[t][[x^{{\rm G}}]]$ по правилу
$$
\delta: \sum_{{\bf m}\in{\mathbb Z}_+^{\kappa}\setminus\{0\}}\alpha_{\bf m}(t)\,x^{m_1r_1+\ldots+m_{\kappa}r_{\kappa}}\mapsto\sum_{{\bf m}\in{\mathbb Z}_+^{\kappa}\setminus\{0\}}\left(m_1r_1+\ldots+m_{\kappa}r_{\kappa}+\frac{d}{dt}\right)\,
\alpha_{\bf m}(t)\,x^{m_1r_1+\ldots+m_{\kappa}r_{\kappa}},
$$
с помощью изоморфизма $\iota: {\mathbb C}[t][[x^{{\rm G}}]]\rightarrow{\mathbb C}[t][[x_1,\ldots,x_{\kappa}]]_*$ естественным образом переносится на ${\mathbb C}[t][[x_1,\ldots,x_{\kappa}]]_*$,
$$
\hat{\delta}: \sum_{{\bf m}\in{\mathbb Z}_+^{\kappa}\setminus\{0\}}\alpha_{\bf m}(t)\,x_1^{m_1}\ldots x_{\kappa}^{m_{\kappa}}\mapsto\sum_{{\bf m}\in{\mathbb Z}_+^{\kappa}\setminus\{0\}}\left(m_1r_1+\ldots+m_{\kappa}r_{\kappa}+\frac{d}{dt}\right)\,\alpha_{\bf m}(t)\,x_1^{m_1}\ldots x_{\kappa}^{m_{\kappa}},
$$
то есть $$\iota\circ\delta=\hat{\delta}\circ\iota.$$

Обозначим
$$
{\bf x}=(x_1,\dots,x_{\kappa}),\qquad {\bf x}^{\bf m}=x_1^{m_1}\ldots x_{\kappa}^{m_{\kappa}}, \qquad {\bf r}=(r_1,\ldots,r_{\kappa}),
$$
 а также
$$
m_1r_1+\ldots+m_\kappa r_\kappa=\langle {\bf m}, {\bf r}\rangle.
$$
Обозначения  ${\bf x}^{\boldsymbol\mu_j},$ ${\bf x}^{\boldsymbol\tau}\;$ и~т.~д. подобны обозначению ${\bf x}^{\bf m}$.

\medskip
Применяя отображение $\iota$ к обеим частям равенства \eqref{reduced} при $u=\psi$  и учитывая свойство (\ref{iotagen}), получим соотношение
\begin{equation}\label{m_reduced}
L\bigl(\lambda_m+\hat{\delta}\bigr)\hat{\psi}+ \widehat{L}\bigl(t,\,{\bf x},\,\lambda_m+\hat{\delta}\bigr)\hat{\psi}+
\widehat{N}\bigl(t,\,{\bf x},\;{\bf x}^{\boldsymbol\tau}\hat{\psi}_0,\ldots,{\bf x}^{\boldsymbol\tau}\hat{\psi}_n\bigr)=0
\end{equation}
для ряда
\begin{equation}\label{m_resh}
\hat{\psi}=\iota(\psi)=\sum\limits_{{\bf m}\in{\mathbb Z}_+^{\kappa}\setminus\{0\}} c_{\bf m}(t)\,{\bf x}^{\bf m} \in\mathbb{C}[t][[x_1,\ldots,x_{\kappa}]]_*,
\end{equation}
где
$\,\hat{\psi}_j=\bigl(\lambda_m+\hat{\delta}\bigr)^j\hat{\psi}$, $\,j=0,\,1,\,\ldots,\,n,\,$ и
$$
\widehat{L}(t,\,{\bf x},\,\zeta)=\sum_{j=0}^{n}\bigl(B_j(t)\,{\bf x}^{{\boldsymbol\mu}_j}+\ldots\bigr)\,\zeta^j\in\mathbb{C}[t]\{x_1,\ldots,x_{\kappa}\}[\zeta],
\qquad {\boldsymbol\mu}_j\in{\mathbb Z}_+^{\kappa}\setminus\{0\}, \quad\mu_j=\langle{\boldsymbol\mu}_j,{\bf r}\rangle,
$$
\begin{eqnarray*}
\widehat{N}\bigl(t,\,{\bf x},\;{\bf x}^{\boldsymbol\tau}\hat{\psi}_0,\ldots,{\bf x}^{\boldsymbol\tau}\hat{\psi}_n\bigr)=
\sum_{\substack{{\bf q}\in\mathbb{Z}_+^{n+1}\\|{\bf q}|\ne1}}\hat{a}_{\bf q}(t,\,{\bf x})\,
\bigl({\bf x}^{\boldsymbol\tau}\hat{\psi}_0\bigr)^{q_0}\cdot\ldots\cdot\bigl({\bf x}^{\boldsymbol\tau}\hat{\psi}_n\bigr)^{q_n}, & &
\hat{a}_{\bf q}\in \mathbb{C}[t]\{x_1,\ldots,x_{\kappa}\} \\ & & \boldsymbol\tau\in{\mathbb Z}_+^{\kappa}\setminus\{0\},\;
\tau=\langle{\boldsymbol\tau},{\bf r}\rangle.
\end{eqnarray*}

\medskip
Итак, уравнение
\begin{eqnarray}\label{m_eq}
L\bigl(\lambda_m+\hat{\delta}\bigr)\hat{u}+ \widehat{L}\bigl(t,\,{\bf x},\,\lambda_m+\hat{\delta}\bigr)\hat{u}+
\widehat{N}\bigl(t,\,{\bf x},\;{\bf x}^{\boldsymbol \tau} \hat{u},\ldots,{\bf x}^{\boldsymbol\tau}(\lambda_m+\hat{\delta})^n\hat{u}\bigr)=0
\end{eqnarray}
имеет  формальное решение \eqref{m_resh}. При этом, в силу однозначности отображения $\iota$ и однозначного определения коэффициентов $c_{k}(t)$ формального  решения \eqref{eq4} уравнения \eqref{reduced}, имеем однозначное определение коэффициентов $c_{\bf m}(t)$  формального решения \eqref{m_resh} уравнения \eqref{m_eq}.

\bigskip
Будем использовать мультииндексные обозначения типа ${\bf m}>{\bf k}$, означающие, что $m_i\geqslant k_i$ при всех $i$ и $|{\bf m}|>|{\bf k}|$.

\bigskip
{\bf Лемма~4.} {\it Каждый многочлен $c_{\bf m}(t)$ в формальном решении \eqref{m_resh} уравнения \eqref{m_eq} имеет степень $$\deg c_{\bf m}\leqslant K\,|{\bf m}|, $$ где $K\geqslant0$ --- некоторая постоянная}.
\medskip

{Д о к а з а т е л ь с т в о.} Напомним, что краткая форма соотношения, которому удовлетворяет ряд $\psi$, имеет вид
$$
x^{-\lambda_m-\nu}F(x,\boldsymbol\varphi_m+x^{\lambda_m}\boldsymbol\psi)=0,
$$
где $\;\boldsymbol\varphi_m=(\varphi_m,\delta\varphi_m,\ldots,\delta^n\varphi_m),\;$ $\boldsymbol\psi=(\psi,(\lambda_m+\delta)\psi,\ldots,(\lambda_m+\delta)^n\psi)$, см. доказательство леммы~2. При этом $F(x,\boldsymbol\varphi_m+x^{\lambda_m}\boldsymbol\psi)$ представляется в виде ряда
$$
\sum_{d\geqslant0}\;x^{d\lambda_m}\;F_d(x,\boldsymbol\varphi_m,\boldsymbol\psi),
$$
$F_d$ --- однородный полином степени $d$ по переменным $\psi_j=(\lambda_m+\delta)^j\psi$, $j=0,1,\ldots,n$:
$$
F_d(x,\boldsymbol\varphi_m,\boldsymbol\psi)=\sum_{\substack{{\bf q}\in{\mathbb Z}_+^{n+1},\\ |{\bf q}|=d}}f_{\bf q} (x,\boldsymbol\varphi_m)\;\psi_0^{q_0}\ldots\psi_n^{q_n},
$$
$f_{\bf q}$ --- голоморфные функции в окрестности точки $0\in{\mathbb C}^{n+2}$, для которых при $|{\bf q}|\ne1$ выполнено соотношение
\begin{equation}\label{d8}
x^{|{\bf q}|\lambda_m-\lambda_m-\nu}f_{\bf q}(x,\boldsymbol\varphi_m)=x^{\tau|{\bf q}|}\,a_{\bf q}(t,x)\in{\mathbb C}[t]\{x^{\rm G}\}.
\end{equation}



Покажем сначала, что для степеней многочленов $\beta_{{\bf q},{\bf i}}(t)$,  $\,{\bf i}\in\mathbb{Z}^{\kappa}_+\setminus\{0\}$, являющихся коэффициентами ряда
${\bf x}^{|{\bf q}|\boldsymbol\tau}\,\hat a_{\bf q}(t,{\bf x})=\sum\limits_{{\bf i}\in\mathbb{Z}_+^{\kappa}\setminus\{0\}}\beta_{{\bf q},\bf i}(t)\,{\bf x}^{\bf i}$ (образа ряда Дюлака $x^{\tau|{\bf q}|}\,a_{\bf q}(t,x)$ под действием отображения $\iota$) справедлива оценка
\begin{equation}\label{d7}
\deg\beta_{{\bf q},{\bf i}}\leqslant K\,|{\bf i}|.
\end{equation}
Обозначим ${\bf i}_0\in{\mathbb Z}_+^\kappa\setminus\{0\}$ мультииндекс, для которого $\lambda_m+\nu=\langle{\bf i}_0,{\bf r}\rangle$.

Пусть
$$
x^{|{\bf q}|\lambda_m}f_{\bf q}(x,\boldsymbol\varphi_m)=\sum_{\substack{p\in{\mathbb Z}_+\setminus\{0\},\\ {\bf s}\in{\mathbb Z}_+^{n+1}\setminus\{0\}}}
\alpha_{{\bf q},p,{\bf s}}\;x^{p+|{\bf q}|\lambda_m}\;\varphi_m^{s_0}(\delta\varphi_m)^{s_1}\ldots(\delta^n\varphi_m)^{s_n}, \qquad
\alpha_{{\bf q},p,{\bf s}}\in{\mathbb C}.
$$
Тогда образ этого ряда Дюлака под действием отображения $\iota$, с одной стороны, имеет вид
$$
\sum_{({\bf p},{\bf s})\in{\mathbb Z}_+^{\kappa+n+1}\setminus\{0\}}
\hat\alpha_{{\bf q},{\bf p},{\bf s}}\;{\bf x}^{\bf p}\;\hat\varphi_m^{s_0}(\hat\delta\hat\varphi_m)^{s_1}\ldots(\hat\delta^n\hat\varphi_m)^{s_n}, \qquad
\hat\alpha_{{\bf q},{\bf p},{\bf s}}\in{\mathbb C}, \quad \hat\varphi_m=\iota(\varphi_m),
$$
а с другой стороны, в силу соотношения (\ref{d8}), --- это ряд $\sum\limits_{{\bf i}\in\mathbb{Z}_+^{\kappa}\setminus\{0\}}\beta_{{\bf q},\bf i}(t)\,{\bf x}^{{\bf i}_0+\bf i}$. Поэтому если
$$
\hat\varphi_m=\sum_{{\bf k}\in{\cal M}}P_{\bf k}(t){\bf x}^{\bf k},
$$
где ${\cal M}\subset{\mathbb Z}_+^\kappa\setminus\{0\}$ --- некоторое конечное множество, то
$\beta_{{\bf q},\bf i}(t)$ --- это линейная комбинация многочленов вида
\begin{equation}\label{eq20}
(P^0_{{\bf k}^0_1}\ldots P^0_{{\bf k}^0_{s_0}})(P^1_{{\bf k}^1_1}\ldots P^1_{{\bf k}^1 _{s_1}})\ldots(P^n_{{\bf k}^n_1}\ldots P^n_{{\bf k}^n_{s_n}}),
\end{equation}
где
$$
P^{j}_{\bf k}(t)=\Bigl(\langle{\bf k},{\bf r}\rangle+\frac{d}{dt}\Bigr)^j P_{\bf k}(t),\qquad j=0,1,\ldots,n,
$$
и мультииндексы ${\bf k}^j_i\in{\cal M}$ удовлетворяют условию ${\bf k}^0_1+{\bf k}^0_2+\ldots+{\bf k}^n_{s_n}\leqslant{\bf i}_0+\bf i$.
Поскольку в \eqref{eq20} все ${\rm deg}\,P^{j}_{\bf k}\leqslant K'=\max_{{\bf k}\in{\cal M}}\deg P_{\bf k}$, то $\deg\beta_{{\bf q},\bf i}\leqslant (s_0+\ldots+s_n)K'$. Но в силу неравенств $s_j\leqslant|{\bf k}^j_1|+|{\bf k}^j_2|+\ldots+|{\bf k}^j_{s_j}|$, $j=0,1,\ldots,n$, получаем, что
$$
\deg\beta_{{\bf q},\bf i}\leqslant K'(|{\bf i}|+|{\bf i}_0|)\leqslant K\,|{\bf i}|, \qquad {\bf i}\in{\mathbb Z}_+^\kappa\setminus\{0\}.
$$

Заметим также, что при $|{\bf q}|=1$ ($q_j=1$, $q_i=0$ при $i\ne j$) выполнено равенство
$$
x^{-\nu}f_{\bf q}(x,\boldsymbol\varphi_m)=A_j+B_j(t)\,x^{\mu_j}+\ldots\in A_j+{\mathbb C}[t]\{x^{\rm G}\},
$$
поэтому оценка, подобная (\ref{d7}), справедлива и для многочленов из коэффициентов оператора $\widehat L(t,{\bf x},\lambda_m+\hat\delta)$.

Перейдём теперь непосредственно к доказательству оценки $\deg c_{\bf m}\leqslant K\,|{\bf m}|$ для коэффициентов ряда \eqref{m_resh}, удовлетворяющего соотношению (\ref{m_reduced}),
$$
L\bigl(\lambda_m+\hat{\delta}\bigr)\hat{\psi}=-\widehat{L}\bigl(t,\,{\bf x},\,\lambda_m+\hat{\delta}\bigr)\hat{\psi}-
\sum_{\substack{{\bf q}\in\mathbb{Z}_+^{n+1}\\|{\bf q}|\ne1}}{\bf x}^{|{\bf q}|\boldsymbol\tau}\hat{a}_{\bf q}(t,\,{\bf x})\,\hat\psi_0^{q_0}\ldots\hat\psi_n^{q_n}.
$$
Будем использовать индукцию по $|{\bf m}|$.

Рассмотрим сначала случай $|{\bf m}|=1$.
Тогда соотношение на коэффициент $c_{\bf m}(t)$ имеет вид
$$
L\Bigl(\lambda_m+\langle{\bf m},{\bf r}\rangle+\frac{d}{dt}\Bigr)c_{\bf m}(t)=\beta_{{\bf 0},\,{\bf m}}(t).
$$
Поскольку $L(\lambda_m+\langle{\bf m},{\bf r}\rangle)\neq0$, то
$$
\deg c_{\bf m}=\deg\beta_{{\bf 0},{\bf m}}\leqslant K.
$$

Пусть теперь $|{\bf m}|\geqslant 2$. В этом случае многочлен $L\bigl(\lambda_m+\langle{\bf m},{\bf r}\rangle+d/dt\bigr)\,c_{\bf m}(t)$ является суммой полиномов вида
\begin{eqnarray}\label{eq21}
-\beta_{{\bf q},\bf s}(t)\bigl(c^0_{{\bf k}^0_1}(t)\ldots c^0_{{\bf k}^0_{q_0}}(t)\bigr)\ldots\bigl(c^n_{{\bf k}^n_1}(t)\ldots c^n_{{\bf k}^n_{q_n}}(t)\bigr),
&\quad& {\bf0}<{\bf s}\leqslant{\bf m}, \\ & & {\bf q}\in\mathbb{Z}_+^{n+1},\;|{\bf q}|\leqslant|{\bf m}|-|{\bf s}|, \nonumber
\end{eqnarray}
где $\deg\beta_{{\bf q},\bf s}\leqslant K\,|{\bf s}|$,
$$
c^{j}_{\bf k}(t)=\Bigl(\lambda_m+\langle{\bf k},{\bf r}\rangle+\frac{d}{dt}\Bigr)^j c_{\bf k}(t),\qquad j=0,1,\ldots,n,
$$
и мультииндексы ${\bf k}^j_i$ удовлетворяют условию
$$
{\bf k}^0_1+{\bf k}^0_2+\ldots+{\bf k}^n_{q_n}={\bf m}-{\bf s}.
$$
Тогда по индукции следует, что в \eqref{eq21} все $\deg c^j_{\bf k}\leqslant K\,|{\bf k}|$, и поэтому суммарная степень \eqref{eq21}, а значит и степень полинома $c_{\bf m}$ не превосходит $K\,|{\bf s}|+K(|{\bf m}|-|{\bf s}|)=K\,|{\bf m}|$. $\hfill\Box$
\medskip

Доказанной леммой мы завершаем обсуждение свойств ряда \eqref{eq2}, удовлетворяющего уравнению \eqref{eq1}.
\medskip

Доказательство теоремы 1 приводится в самом конце работы. Оно опирается на теорему о неявном отображении для банаховых пространств, основные идеи  применения которой заимствованы из работы \cite{Mal} и стали уже стандартными.
Докажем сначала некоторые вспомогательные утверждения, необходимые для доказательства теоремы~1.
\medskip

Во множестве мультииндексов ${\bf m}\in{\mathbb Z}_+^\kappa\setminus\{0\}$, удовлетворяющих условию ${\rm Re}\,\langle{\bf m},{\bf r}\rangle>\tau$, имеется ненулевое число минимальных (относительно отношения порядка ''$>$'', определённого перед леммой 4) элементов.\footnote{Утверждение, известное под названием леммы Диккенса, см., например, \cite{Dick}.} Обозначим их $\;{\bf m}^{(1)},\ldots,{\bf m}^{(d)}\;$ и рассмотрим следующие нормированные пространства $H^j$, $j=0,1,\ldots,\ell,$ формальных рядов из $\mathbb{C}[t][[x_1,\ldots,x_{\kappa}]]_*$:
$$
H^j=\Bigl\{\;\eta=\sum\limits_{{\bf m}\in\mathbb{Z}_+^\kappa\setminus\{0\}}C_{\bf m}(t)\,{\bf x}^{\bf m}\quad\mid\quad \deg C_{\bf m}\leqslant \mathcal{K}\,
|{\bf m}|,\quad \|\eta\|_j<+\infty\Bigr\},
$$
где $\;\mathcal{K}=2\,K\,\max\limits_{i=1,\ldots,d}|{\bf m}^{(i)}|,\;$ ${\bf m}^{(i)}\,$ определены выше, $\;K$ --- в лемме 4,
$$
\|\eta\|_j=\sum\limits_{{\bf m}\in\mathbb{Z}_+^\kappa\setminus\{0\}}\frac{(|\lambda_m+\langle{\bf m},{\bf r}\rangle|+\mathcal{K}\,|{\bf m}|)^j}
{|\Gamma\left(\langle{\bf m},{\bf r}\rangle/s\right)|}\,\|C_{\bf m}\|_R.
$$
Поскольку пространства многочленов фиксированной степени (как конечномерные) полны по норме $\|\cdot\|_R$, то пространства $H^j$ являются банаховыми (их полнота проверяется так же, как полнота пространства $l_1$ абсолютно суммируемых числовых последовательностей). Также пространства $H^0,\ldots, H^\ell$ обладают следующими свойствами:
\begin{itemize}
\item[] $\qquad H^\ell\subset H^{\ell-1}\subset\ldots\subset H^0,$
\item[] $\qquad \lambda_m+\hat{\delta}: H^j\rightarrow H^{j-1}$ --- непрерывный линейный оператор, $j=1,\ldots,\ell.$
\end{itemize}
В частности,  оператор $\;(\lambda_m+\hat{\delta})^{\ell}\;$  действует из $H^{\ell}$ в $H^0.$ Но при $j>\ell$ оператор $(\lambda_m+\hat{\delta})^{j}$ может отобразить $H^{\ell}$ вне $H^0$. Это можно исправить, рассмотрев такой оператор, но с указанным в следующей лемме множителем.

\bigskip
{\bf Лемма 5.} {\it Пусть  $\;{\bf l}\in{\mathbb Z}_+^\kappa\setminus\{0\}\;$ такой, что $${\rm Re}\,\langle{\bf l},{\bf r}\rangle\geqslant (j-\ell)s,\qquad j>\ell.$$ Тогда,  если $a\in{\mathbb C}[t]$, $\deg a\leqslant\mathcal{K}\,|{\bf l}|$, то линейный оператор
$$
a(t)\,{\bf x}^{\bf l}\,(\lambda_m+\hat{\delta})^j: H^\ell\rightarrow H^0
$$
является непрерывным.}
\bigskip

{Д\,о\,к\,а\,з\,а\,т\,е\,л\,ь\,с\,т\,в\,о.} Рассмотрим ряд
$$
\;\eta=\sum\limits_{{\bf m}\in\mathbb{Z}_+^\kappa\setminus\{0\}}C_{\bf m}(t)\,{\bf x}^{\bf m} \in H^\ell.
$$
Обозначим
$$
\mathcal{C}_{\bf m}(t)=\frac{(|\lambda_m+\langle{\bf m},{\bf r}\rangle|+\mathcal{K}\,|{\bf m}|)^{\ell}}{|\Gamma\left(\langle{\bf m},{\bf r}\rangle/s\right)|}\,
C_{\bf m}(t),
$$
тогда числовой ряд $\;\sum\limits_{{\bf m}\in\mathbb{Z}_+^\kappa\setminus\{0\}}\|\mathcal{C}_{\bf m}\|_R\;$ сходится, что следует из определения $H^{\ell}$.

\bigskip
Также рассмотрим ряд
$$
a(t)\,{\bf x}^{\bf l}\,(\lambda_m+\hat{\delta})^j\,\eta=\sum_{{\bf m}\in\mathbb{Z}_+^\kappa\setminus\{0\}}a(t)\left(\lambda_m+\langle{\bf m},{\bf r}\rangle+ \frac{d}{dt}\right)^j C_{\bf m}(t)\;{\bf x}^{{\bf m}+{\bf l}}=
$$
$$
=\sum_{\substack{{\bf m}\in\mathbb{Z}_+^\kappa\setminus\{0\}\\{\bf m}>{\bf l}}}a(t)\left(\lambda_m+\langle{\bf m}-{\bf l},{\bf r}\rangle+ \frac{d}{dt}\right)^j\,C_{{\bf m}-{\bf l}}(t)\,{\bf x}^{\bf m}=
$$
$$
=\sum_{\substack{{\bf m}\in\mathbb{Z}_+^\kappa\setminus\{0\}\\{\bf m}>{\bf l}}}a(t)\left(\lambda_m+\langle{\bf m}-{\bf l},{\bf r}\rangle+ \frac{d}{dt}\right)^j\,\frac{|\Gamma(\langle{\bf m}-{\bf l},{\bf r}\rangle/s)|}{(|\lambda_m+\langle{\bf m}-{\bf l},{\bf r}\rangle|+\mathcal{K}\,|{\bf m}-{\bf l}|) ^{\ell}}\; \mathcal{C}_{{\bf m}-{\bf l}}(t)\;{\bf x}^{\bf m}.
$$

Покажем, что $a(t)\,{\bf x}^{\bf l}(\lambda_m+\hat{\delta})^j\eta\in H^0$. Требуемая оценка на степени коэффициентов этого ряда выполнена, поскольку $\,\deg\mathcal{C}_{\bf m}=\deg C_{\bf m}\leqslant\mathcal{K}\,|{\bf m}|$ согласно определению пространств $H^j$, а $\deg a\leqslant  \mathcal{K}\,|{\bf l}|$ по условию леммы. Остаётся доказать сходимость числового ряда
\begin{eqnarray}\label{h0}
\sum_{\substack{{\bf m}\in\mathbb{Z}_+^\kappa\setminus\{0\}\\{\bf m}>{\bf l}}}\frac{|\Gamma(\langle {\bf m}-{\bf l},{\bf r}\rangle/s)|}{|\Gamma(\langle{\bf m},{\bf r}\rangle/s)|}\;\;\frac{\displaystyle\left\|a(t)\,\Bigl(\lambda_m+\langle{\bf m}-{\bf l},{\bf r}\rangle+\frac{d}{dt}\Bigr)^j\mathcal{C}_{{\bf m}-{\bf l}}\right\|_R} {(|\lambda_m+\langle{\bf m}-{\bf l},{\bf r}\rangle|+\mathcal{K}\,|{\bf m}-{\bf l}|)^\ell}.
\end{eqnarray}

Во-первых, с учётом оценки
$$
\Bigl\|\left(\frac d{dt}\right)^i\mathcal{C}_{{\bf m}-\bf l}\Bigr\|_R\leqslant(\deg\mathcal{C}_{{\bf m}-\bf l})^i\,\|\mathcal{C}_{{\bf m}-\bf l}\|_R\leqslant
(\mathcal{K}\,|{\bf m}-{\bf l}|)^i\,\|\mathcal{C}_{{\bf m} -\bf l}\|_R
$$
получаем
$$
\left\|\left(\lambda_m+\langle{\bf m}-{\bf l},{\bf r}\rangle+\frac{d}{dt}\right)^{j}\mathcal{C}_{{\bf m}-{\bf l}}\right\|_R\leqslant \left(|\lambda_m+\langle{\bf m}- {\bf l},{\bf r}\rangle|+\mathcal{K}\,|{\bf m}-{\bf l}|\right)^j \|\mathcal{C}_{{\bf m}-{\bf l}}\|_R.
$$
Во-вторых, из асимптотических свойств гамма-функции (см., например, п. 1.18 в \cite{Batm}) следует, что $\,\Gamma(z)/\Gamma(z+a)\sim z^{-a}\,$ при ${\rm Re}\,z\rightarrow\infty$ и фиксированном ${\rm Re}\, a>0$. Поэтому, в силу отделимости ${\rm Re}\langle {\bf m}-{\bf l},{\bf r}\rangle$ от нуля, имеем
$$
\left|\frac{\Gamma(\langle {\bf m}-{\bf l},{\bf r}\rangle/s)}{\Gamma(\langle{\bf m},{\bf r}\rangle/s)}\right| \,\leqslant\, \widetilde A_1\,|\langle {\bf m}-
{\bf l},{\bf r}\rangle^{-\langle{\bf l},{\bf r}\rangle/s}|,
$$
для некоторого $\widetilde A_1>0$ и всех ${\bf m}>{\bf l}$, где ${\bf l}$ фиксировано.  Отметим также выполнение неравенств
$$
|\lambda_m+\langle{\bf m}-{\bf l},{\bf r}\rangle|\leqslant \alpha\,|\langle{\bf m}- {\bf l},{\bf r}\rangle|, \qquad |\langle{\bf m}- {\bf l},{\bf r}\rangle| \geqslant{\rm Re}\,\langle{\bf m}-{\bf l},{\bf r}\rangle\geqslant\beta\,|{\bf m}-{\bf l}|
$$
при подходящих $\alpha,\beta>0$.

Таким образом, с учётом сделанных замечаний и оценки $\|P\,Q\|_R\leqslant\|P\|_R\cdot\|Q\|_R$ ($P,$ $Q$ --- произвольные многочлены), ряд
\begin{equation}\label{eq22}
\widetilde A_2\sum\limits_{\substack{{\bf m}\in\mathbb{Z}_+^\kappa\setminus\{0\}\\{\bf m}>{\bf l}}}\;\bigl|\langle{\bf m}-{\bf l},{\bf r}\rangle^{j-\ell-\langle{\bf l},{\bf r}\rangle/s}\bigr|\,\|\mathcal{C}_{{\bf m}-{\bf l}}\|_R
\end{equation}
является мажорантным для ряда (\ref{h0}) при некотором подходящем $\widetilde A_2>0$. Ряд \eqref{eq22} сходится, поскольку при $|{\bf m}|\rightarrow \infty$  выражение $|\langle {\bf m}-{\bf l},{\bf r}\rangle|\rightarrow\infty$, а $j-\ell-{\rm Re}(\langle{\bf l},{\bf r}\rangle/s)\leqslant0$ в предположениях леммы. Отсюда следует сходимость ряда \eqref{h0} и непрерывность линейного оператора $a(t)\,{\bf x}^{\bf l}(\lambda_m+\hat{\delta})^j$. Действительно,
$$
\left\|a(t)\,{\bf x}^{\bf l}(\lambda_m+\hat{\delta})^j\eta\right\|_0\leqslant \widetilde A_2\sum\limits_{\substack{{\bf m}\in\mathbb{Z}_+^\kappa\setminus\{0\}\\
{\bf m}>{\bf l}}}\; \left|\langle {\bf m}-{\bf l},{\bf r}\rangle^{j-\ell-\langle{\bf l},{\bf r}\rangle/s}\right|\,\|\mathcal{C}_{{\bf m}-{\bf l}}\|_R\leqslant
$$
$$
\leqslant \widetilde A_3\sum_{{\bf m}\in\mathbb{Z}_+^\kappa\setminus\{0\}}\|\mathcal{C}_{\bf m}\|_R=\widetilde A_3
\sum_{{\bf m}\in\mathbb{Z}_+^\kappa\setminus\{0\}} \frac{(|\lambda_m+\langle{\bf m},{\bf r}\rangle|+\mathcal{K}\,|{\bf m}|)^{\ell}}{|\Gamma\left(\langle{\bf m},
{\bf r}\rangle/s\right)|}\,\|C_{\bf m}\|_R=\widetilde A_3\;\|\eta\|_\ell, \qquad \widetilde A_3>0.
$$
{\hfill $\Box$}
\medskip

{\bf Лемма 6.} {\it Для произвольных $\eta_1$, $\eta_2\in H^0$ произведение $\eta_1\cdot\eta_2\in H^0$, при этом
$$
\|\eta_1\eta_2\|_0\leqslant C\|\eta_1\|_0\|\eta_2\|_0, \qquad  C\geqslant1.
$$}

\bigskip
{Д\,о\,к\,а\,з\,а\,т\,е\,л\,ь\,с\,т\,в\,о.}  Пусть
$$
\eta_1=\sum\limits_{{\bf m}\in\mathbb{Z}_+^\kappa\setminus\{0\}}a_{\bf m}(t)\,{\bf x}^{\bf m}, \qquad
\eta_2=\sum\limits_{{\bf m}\in\mathbb{Z}_+^\kappa\setminus\{0\}}b_{\bf m}(t)\,{\bf x}^{\bf m}.
$$
Тогда
$$
\eta_1\eta_2=\sum_{{\bf m}\in\mathbb{Z}_+^\kappa\setminus\{0\}}\biggl(\sum_{{\bf 0}<{\bf i}<{\bf m}}\,a_{\bf i}(t)\,b_{{\bf m}-{\bf i}}(t)\biggr)\,{\bf x}^{\bf m},
$$
$$
\|\eta_1\eta_2\|_0=\sum_{{\bf m}\in\mathbb{Z}_+^\kappa\setminus\{0\}}\frac{\displaystyle\biggl\|\sum_{{\bf 0}<{\bf i}<{\bf m}}a_{\bf i}\,b_{{\bf m}-{\bf i}} \biggr\|_R}{\displaystyle|\Gamma(\langle{\bf m},{\bf r}\rangle/s)|}.
$$
Пользуясь известным соотношением между $\Gamma$-функцией и $B$-функцией,
$$
\frac{\Gamma(\langle{\bf i},{\bf r}\rangle/s)\Gamma(\langle{\bf m}-{\bf i},{\bf r}\rangle/s)}{\Gamma(\langle{\bf m},{\bf r}\rangle/s)}=\int_0^1 \displaystyle(1-t)^{\langle{\bf i},{\bf r}\rangle/s-1}\,t^{\langle{\bf m}-{\bf i},{\bf r}\rangle/s-1}\,dt,
$$
в силу отделимости ${\rm Re}\,\langle{\bf i},{\bf r}\rangle$ и ${\rm Re}\,\langle{\bf m}-{\bf i},{\bf r}\rangle$ от нуля получаем
$$
\left|\frac{\Gamma(\langle{\bf i},{\bf r}\rangle/s)\Gamma(\langle{\bf m}-{\bf i},{\bf r}\rangle/s)}{\Gamma(\langle{\bf m},{\bf r}\rangle/s)}\right|\leqslant
\int_0^1\displaystyle(1-t)^{{\rm Re}\,\langle{\bf i},{\bf r}\rangle/s-1}\,t^{{\rm Re}\,\langle{\bf m}-{\bf i},{\bf r}\rangle/s-1}\,dt\leqslant C
$$
для любого ${\bf 0}<{\bf i}<{\bf m}$, при некотором $C\geqslant1$. Тогда из неравенства
$$
\frac1{|\Gamma(\langle{\bf m},{\bf r}\rangle/s)|}\leqslant\frac C{|\Gamma(\langle{\bf i},{\bf r}\rangle/s)\,\Gamma(\langle{\bf m}-{\bf i},{\bf r}\rangle/s)|}
$$
следует
$$
\|\eta_1\eta_2\|_0\leqslant C\sum_{{\bf m}\in\mathbb{Z}_+^\kappa\setminus\{0\}}\sum_{{\bf 0}<{\bf i}<{\bf m}}\frac{\|a_{\bf i}\|_R}
{|\Gamma(\langle{\bf i},{\bf r}\rangle/s)|}\,\frac{\|b_{{\bf m}-{\bf i}}\|_R}{|\Gamma(\langle{\bf m}-{\bf i},{\bf r}\rangle/s)|}= C\|\eta_1\|_0\|\eta_2\|_0.
$$
Остаётся заметить, что требуемая оценка на степени коэффициентов ряда $\eta_1\eta_2$ также выполнена, поскольку
$\deg a_{\bf i}\,b_{{\bf m}-{\bf i}}\leqslant \mathcal{K}\,|{\bf m}|$. {\hfill $\Box$}

\bigskip
Обозначим
$$
{\boldsymbol \psi}=(\psi_0,\ldots,\psi_n), \quad {\bf s}=(\|\psi_0\|_0,\ldots,\|\psi_n\|_0),\qquad {\bf p}=(p_1,\ldots,p_{\kappa})\in\mathbb{Z}^{\kappa}_+\setminus\{0\}.
$$
Обозначения ${\bf u}^{\bf q}$, ${\boldsymbol \psi}^{\bf q}$, ${\bf s}^{\bf q}$ подобны обозначению ${\bf x}^{\bf m}$.

\bigskip
{\bf Лемма 7.} {\it Пусть $F(t,{\bf x},{\bf u})$ --- полином по $t$ и голоморфная  по переменным $\bf x$, $\bf u$ в окрестности точки $(0,0)\in\mathbb{C}^{\kappa+n+1}$ функция,  представимая рядом
$$
F(t,{\bf x},{\bf u})=\sum\limits_{({\bf p},{\bf q})\in\mathbb{Z}^{\kappa+n+1}_+\setminus\{0\}}a_{{\bf p},{\bf q}}(t)\,{\bf x}^{\bf p}{\bf u}^{\bf q},
$$
где $a_{{\bf p},{\bf q}}\in{\mathbb C}[t]$ и $\deg a_{{\bf p},{\bf q}}\leqslant\mathcal{K}\,|{\bf p}|$. Пусть также
$\psi_0,$ $\ldots,$ $\psi_n\in H^{0}$, $\rho>0$. Тогда для любого $\varepsilon>0$ будет выполнено $F(t,{\rho}\,{\bf x},{\boldsymbol\psi})\in H^0$ и $\left\|F(t,{\rho}\,{\bf x},{\boldsymbol\psi})\right\|_0<\varepsilon$ при достаточно малых $\|\psi_j\|_0$, $\rho$.}

\medskip
{Д о к а з а т е л ь с т в о.} Рассмотрим последовательность $\{F_1,\,F_2,\ldots\}\subset H^0$:
$$
F_N=\sum\limits_{\substack{{({\bf p},{\bf q})\in\mathbb{Z}^{\kappa+n+1}_+\setminus\{0\}}\\{|{\bf p}|+|{\bf q}|\leqslant N}}}a_{{\bf p},{\bf q}}(t)\,{\rho^{|{\bf p}|}}\,{\bf x}^{\bf p} {\boldsymbol\psi}^{\bf q}
$$
 ($F_N\in H^0$ согласно лемме 6). Покажем, что она является фундаментальной в $H^0$ ({при достаточно малых $\|\psi_j\|_0$, $\rho>0$}). 
При всех $N_1<N_2$ ввиду леммы 6 имеем
$$
\|F_{N_2}-F_{N_1}\|_0=\Bigl\|\sum\limits_{\substack{{({\bf p},{\bf q})\in\mathbb{Z}^{\kappa+n+1}_+\setminus\{0\}}\\{N_1<|{\bf p}|+|{\bf q}|\leqslant N_2}}}
a_{{\bf p},{\bf q}}(t)\;{\rho^{|{\bf p}|}}\;{\bf x}^{\bf p}\;{\boldsymbol\psi}^{\bf q}\Bigr\|_0\leqslant$$
$$\leqslant\sum\limits_{\substack{{({\bf p},{\bf q}) \in\mathbb{Z}^{\kappa+n+1}_+\setminus\{0\}}\\{N_1<|{\bf p}|+|{\bf q}|\leqslant N_2}}}\frac{\|a_{{\bf p},{\bf q}}\|_R}{\bigl|\Gamma(\langle{\bf p},{\bf r}\rangle/s)\bigr|}\;\;{\rho^{|{\bf p}|}\;C^{|{\bf q}|}}\;{{\bf s}^{\bf q}}\rightarrow 0
$$
при $N_1,$ $N_2\rightarrow +\infty$, поскольку числовой ряд $\;\displaystyle\sum\limits_{({\bf p},{\bf q})\in\mathbb{Z}^{\kappa+n+1}_+\setminus\{0\}}
\frac{\|a_{{\bf p},{\bf q}}\|_R}{\bigl|\Gamma(\langle{\bf p},{\bf r}\rangle/s)\bigr|}\;\;\rho^{|{\bf p}|}\;C^{|{\bf q}|}\;{\bf s}^{\bf q}\;$ сходится при достаточно малых $\|\psi_j\|_0$, $\rho>0$. Таким образом, последовательность $\{F_1,\,F_2,\ldots\}$ имеет предел в $H^0$. Покажем теперь, что $\,F(t,{\rho}\,{\bf x},{\boldsymbol\psi})\in H^0\,$ и $\,\|F(t,{\rho}\,{\bf x}, {\boldsymbol\psi})\|_0$ мало при малых $\,\|\psi_j\|_0$, $\rho>0$.

Запишем
$$
F(t,{\rho}\,{\bf x},{\boldsymbol \psi})=\sum_{{\bf p}\in\mathbb{Z}_+^{\kappa}\setminus\{0\}}\alpha_{\bf p}(t)\,{\bf x}^{\bf p}
$$
и заметим, что
$$
F(t,{\rho}\,{\bf x},{\boldsymbol \psi})-F_N=\sum\limits_{\substack{{({\bf p},{\bf q})\in\mathbb{Z}^{\kappa+n+1}_+\setminus\{0\}}\\{|{\bf p}|+|{\bf q}|>N}}}
a_{{\bf p},{\bf q}}(t)\,\rho^{|{\bf p}|}\,{\bf x}^{\bf p}{\boldsymbol\psi}^{\bf q}=
\sum_{\substack{{{\bf p}\in\mathbb{Z}_+^{\kappa}\setminus\{0\}}\\{|{\bf p}|>N}}}\beta^{(N)}_{\bf p}(t)\,{\bf x}^{\bf p},
$$
то есть
$$
F_N=\sum_{\substack{{{\bf p}\in\mathbb{Z}_+^{\kappa}\setminus\{0\}}\\{|{\bf p}|\leqslant N}}}\alpha_{\bf p}(t)\,{\bf x}^{\bf p}+\sum_{\substack{{{\bf p}\in \mathbb{Z}_+^{\kappa}\setminus\{0\}}\\{|{\bf p}|> N}}}\tilde{\alpha}^{(N)}_{\bf p}(t)\,{\bf x}^{\bf p}.
$$
Поскольку все $F_N\in H^0$, то $\deg\alpha_{\bf p}\leqslant\mathcal{K}\,|{\bf p}|$. Также имеем
$$\|F_N\|_0=\sum_{\substack{{{\bf p}\in\mathbb{Z}_+^{\kappa}\setminus\{0\}}\\{|{\bf p}|\leqslant N}}}\frac{\|\alpha_{\bf p}\|_R}{|\Gamma(\langle{\bf p},{\bf r} \rangle/s)|}+\sum_{\substack{{{\bf p}\in\mathbb{Z}_+^{\kappa}\setminus\{0\}}\\{|{\bf p}|>N}}}\frac{\|\tilde{\alpha}^{(N)}_{\bf p}\|_R}{|\Gamma(\langle{\bf p},{\bf r} \rangle/s)|},
$$
поэтому
$$
\sum_{\substack{{{\bf p}\in\mathbb{Z}_+^{\kappa}\setminus\{0\}}\\{|{\bf p}|\leqslant N}}}\frac{\|\alpha_{\bf p}\|_R}{|\Gamma(\langle{\bf p},{\bf r}\rangle/s)|} \leqslant\|F_N\|_0.
$$
Следовательно,
$$
\|F(t,{\rho}\,{\bf x},{\boldsymbol\psi})\|_0=\lim_{N\rightarrow+\infty}\sum_{\substack{{{\bf p}\in\mathbb{Z}_+^{\kappa}\setminus\{0\}}\\
{|{\bf p}|\leqslant N}}}\frac{\|\alpha_{\bf p}\|_R}{|\Gamma(\langle{\bf p},{\bf r}\rangle/s)|}\leqslant
$$
$$
\leqslant
\lim_{N\rightarrow+\infty}\|F_N\|_0\leqslant\sum\limits_{({\bf p},{\bf q})\in\mathbb{Z}^{\kappa+n+1}_+\setminus\{0\}}\frac{\|a_{{\bf p},{\bf q}}\|_R} {\bigl|\Gamma(\langle{\bf p},{\bf r}\rangle/s)\bigr|} \rho^{|{\bf p}|}C^{|{\bf q}|}\,{\bf s}^{\bf q},
$$
а последняя величина мала при малых $\|\psi_j\|_0$, $\rho>0$  (равна нулю при $\rho=0, {\bf s}=0$). $\hfill\Box$

\bigskip
Далее мы будем использовать следующую теорему.

\medskip
{\bf Теорема} {\rm (о неявном отображении для банаховых пространств\footnote{Здесь приводится ослабленный вариант теоремы, поскольку нас интересует только {\it существование} неявного отображения.}, см., например, \cite[гл. X, \S2, п. 1]{KolmFomin}).} {\it Пусть $\cal E$, $\cal F$, $\cal H$ --- банаховы пространства, $V$ --- открытое подмножество прямого произведения ${\cal E}\times{\cal F}$, и $h: V\rightarrow\cal H$ --- непрерывное в точке $(x_0, y_0)\in V$ отображение такое, что
\begin{itemize}
\item[] $h(x_0,y_0)=0$;
\item[] производная $h'_y(x,y)$ $($непрерывное линейное отображение из $\cal F$ в $\cal H)$ определена при всех $(x,y)\in V$ и непрерывна в $(x_0,y_0)$;
\item[] $h'_y(x_0,y_0)$ --- линейный гомеоморфизм из $\cal F$ в $\cal H$.
\end{itemize}
Тогда имеется окрестность $U_0\subset\cal E$ точки $x_0$ и отображение $g: U_0\rightarrow\cal F$ такие, что $g(x_0)=y_0$, $(x,g(x))\in V$, и $h(x,g(x))=0$
для любого $x\in U_0$.}

\bigskip
{Д\,о\,к\,а\,з\,а\,т\,е\,л\,ь\,с\,т\,в\,о\,\, т\,е\,о\,р\,е\,м\,ы\, 1.\, }
Применим данную теорему к банаховым пространствам $\mathbb C$, $H^\ell$, $H^0$ и отображению $h:V\subset{\mathbb C}\times H^\ell\rightarrow H^0$, определённому в окрестности $V$ точки $(0,0)\in{\mathbb C}\times H^\ell$ следующим образом:
$$
h: (\rho,\eta)\mapsto L(\lambda_m+\hat{\delta})\eta+\widehat{L}(t,\rho\,{\bf x},\lambda_m+\hat{\delta})\eta+
\widehat N(t,\rho\,{\bf x},\,\rho^{|{\boldsymbol\tau}|}{\bf x}^{\boldsymbol \tau}\,\eta,\ldots,
\rho^{|{\boldsymbol \tau}|}{\bf x}^{\boldsymbol\tau}\,(\lambda_{m}+\hat{\delta})^n\eta),
$$
где $L$, $\widehat L$ и $\widehat N$ взяты из соотношения (\ref{m_reduced}). Проверим правомерность её применения.

\medskip
\begin{itemize}
\item[$\circ$] {Сначала покажем, что $\widehat{L}(t, \rho\, {\bf x},\lambda_m+\hat{\delta})\eta\in H^0$.
    Представим ряды
    $$
    B_j(t){\bf x}^{\boldsymbol\mu_j}+\ldots\in{\mathbb C}[t]\{x_1,\ldots,x^\kappa\},
    $$
    являющиеся коэффициентами линейного дифференциального оператора $\widehat{L}\bigl(t,\,{\bf x},\,\lambda_m+\hat\delta\bigr)$, в виде
    $$
    \sum_{{\rm Re}\,\mu_j\leqslant{\rm Re}\,\langle{\bf m},{\bf r}\rangle\leqslant\tau}\beta_{j,\bf m}(t){\bf x}^{\bf m}+
    {\bf x}^{{\bf m}^{(1)}}\sum_{{\bf m}\in{\mathbb Z}_+^\kappa\setminus\{0\}}\beta^{(1)}_{j,\bf m}(t){\bf x}^{\bf m}+\ldots+
    {\bf x}^{{\bf m}^{(d)}}\sum_{{\bf m}\in{\mathbb Z}_+^\kappa\setminus\{0\}}\beta^{(d)}_{j,\bf m}(t){\bf x}^{\bf m}.
    $$

    Как следует из доказательства леммы 4,
    $$
    \deg\beta_{j,\bf m}\leqslant K\,|{\bf m}|\leqslant\mathcal{K}\,|{\bf m}|, \quad
    \deg\beta^{(i)}_{j,\bf m}\leqslant K\,|{\bf m}+{\bf m}^{(i)}|\leqslant2K\,|{\bf m}|\,|{\bf m}^{(i)}|\leqslant\mathcal{K}\,|{\bf m}|,
    $$
    то есть все $\displaystyle\sum\limits_{{\bf m}\in{\mathbb Z}_+^\kappa\setminus\{0\}}\beta^{(i)}_{j,\bf m}(t){\bf x}^{\bf m}\in H^0$. Кроме того, согласно лемме 5 все $\beta_{j,\bf m}(t){\bf x}^{\bf m}\,(\lambda_m+\hat{\delta})^j\eta\in H^0$, так как по замечанию~1 $\;{\rm Re}\,\langle{\bf m},{\bf r}\rangle \geqslant{\rm Re}\,\mu_j\geqslant (j-\ell)s$ при $j>\ell$. Также согласно лемме~5 все $\,{\bf x}^{{\bf m}^{(i)}}\,(\lambda_m+\hat{\delta})^j\eta\in H^0$, поскольку ${\rm Re}\,\langle{\bf m}^{(i)},{\bf r}\rangle>\tau=(n-\ell)s$. Таким образом, согласно леммам 6, 7 отображение $h(\rho,\eta)$ корректно определено при малых $\rho$ и $\|\eta\|_\ell$.}

\item[$\circ$] Отображение $h$ непрерывно в точке $(0,0)\in{\mathbb C}\times H^\ell$ в силу лемм 5, 6 и 7, при этом $h(0,0)=0$.

\item[$\circ$] Непосредственно проверяется (с использованием формальной формулы Тейлора и лемм 5, 6 и 7), что производная $h'_\eta(\rho,\eta)$ определена в $V$ соответствием

\medskip\noindent
$
h'_\eta: (\rho,\eta)\mapsto L(\lambda_m+\hat{\delta})+\widehat{L}(t, \rho\, {\bf x},\lambda_m+\hat{\delta})+\hfill\phantom{a}
$
$$\eqno\sum_{j=0}^n\widehat N'_{u_j}(t,\,\rho\,{\bf x},\,\rho^{|\boldsymbol\tau|}{\bf x}^{\boldsymbol\tau}\,\eta,\ldots,
\rho^{|\boldsymbol\tau|}{\bf x}^{\boldsymbol\tau}\,(\lambda_m+\hat{\delta})^n\eta)\;\rho^{|\boldsymbol\tau|}{\bf x}^{\boldsymbol\tau}(\lambda_m+\hat{\delta})^j,
$$
в частности, $h'_\eta(0,0)=L(\lambda_m+\hat{\delta})$ и разность $h'_\eta(\rho,\eta)-h'_\eta(0,0)$ мала при малых  $\rho$ и $\|\eta\|_\ell$ по норме пространства ${\cal L}(H^\ell,H^0)$ непрерывных линейных отображений $H^\ell\rightarrow H^0$. Поэтому производная $h'_\eta$ непрерывна в точке
$(0,0)\in{\mathbb C}\times H^\ell$.

\item[$\circ$] Наконец, $h'_\eta(0,0)=L(\lambda_m+\hat{\delta})$ --- непрерывное линейное отображение из  $H^\ell$ в $H^0$ (непрерывность операторов $(\lambda_m+\hat{\delta})^j$ обсуждалась ранее), обладающее непрерывным обратным. Проверим последнее свойство. Сначала сделаем некоторые вспомогательные оценки.

    Пусть $\displaystyle{\Delta=\lambda_m+\langle{\bf m},{\bf r}\rangle+d/dt}$. Заметим, что для любого многочлена $C(t)$ степени, не превышающей $\mathcal{K}\,|{\bf m}|$, имеет место неравенство
      \begin{equation}\label{th_eq1}
      \left\|\frac{d}{dt}C(t)\right\|_R\leqslant \frac{\mathcal{K}|{\bf m}|}{R}\|C(t)\|_R.
\end{equation}
Так как $\displaystyle C(t)=\frac{1}{{\lambda_m+}\langle{\bf m},{\bf r}\rangle}\Bigl({-}\frac{d}{dt}+\Delta\Bigr)C(t)$, то с учётом \eqref{th_eq1} получаем неравенство
$$
\|C(t)\|_R\leqslant \frac{\mathcal{K}}{R}\,\frac{|{\bf m}|}{|{\lambda_m+}\langle{\bf m},{\bf r}\rangle|}\,\|C(t)\|_R+ \frac{1}{|{\lambda_m+}\langle{\bf m},{\bf r}\rangle|} \,\|\Delta C(t)\|_R.
$$
Обозначим через $\displaystyle\;\theta={\sup_{\bf m\in\mathbb{Z}_+^{\kappa}\setminus\{0\}}} {|{\bf m}|}/{|{\lambda_m+}\langle{\bf m},
{\bf r}\rangle|}\,$ и $\displaystyle\,\Theta=\theta/\Bigl(1-\frac{\mathcal{K}}R\theta\Bigr).$ Тогда для $R>\mathcal{K}\theta$ {(этим условием и определяется выбор величины $R$)} справедливо неравенство
\begin{equation}\label{th_eq2}
    \|C(t)\|_R\leqslant \Theta \frac{\|\Delta C(t)\|_R}{|{\bf m}|}\leqslant{\ldots\leqslant}\Theta^j\frac{\|\Delta^j C(t)\|_R}{|{\bf m}|^j},\qquad j\geqslant 1.
\end{equation}

Поскольку
$$
\Delta^{\ell}C(t)=A_{\ell}^{-1}\left(L\Bigl(\lambda_m+\langle{\bf m},{\bf r}\rangle+\frac{d}{dt}\Bigr)C(t)-A_{\ell-1}\Delta^{\ell-1}C(t)-\ldots-A_0C(t)\right),
$$
то с учётом \eqref{th_eq2} получаем цепочку неравенств
$$
\|\Delta^{\ell}C(t)\|_R\leqslant \alpha_\ell\left\|L\left(\lambda_m+\langle{\bf m},{\bf r}\rangle+d/dt\right)C(t)\right\|_R+ \alpha_{\ell-1}\|\Delta^{\ell-1}C(t)\|_R+\ldots+\alpha_0\|C(t)\|_R\leqslant
$$
$$
\leqslant\alpha_\ell \left\|L\left(\lambda_m+\langle{\bf m},{\bf r}\rangle+d/dt\right)C(t)\right\|_R+
\alpha_{\ell-1}\,\Theta\,\frac{\|\Delta^{\ell}C(t)\|_R}{|{\bf m}|}+\ldots+\alpha_0\,\Theta^{\ell}\,\frac{\|\Delta^{\ell}C(t)\|_R}{|{\bf m}|^{\ell}},
$$
где $\alpha_{\ell}=|A_{\ell}^{-1}|$ и $\alpha_j=|A_{\ell}^{-1}A_{j}|$, $j=1,\ldots,\ell-1$. Отсюда следует, что
$$
\left(1-\frac{\alpha_{\ell-1}{\Theta}}{|{\bf m}|}-\ldots-\frac{\alpha_{0}{\Theta^\ell}}{|{\bf m}|^{\ell}}\right) \|\Delta^{\ell}C(t)\|_R\leqslant\alpha_{\ell}\left\|L(\lambda_m+\langle{\bf m},{\bf r}\rangle+d/dt)C(t)\right\|_R.
$$
Существует $\mu>0$ такое, что при  всех $|{\bf m}|>\mu$ величина, стоящая в скобках, положительна и отделена от нуля. Отсюда следует, что имеется такое $\alpha>0$,  не зависящее от  ${\bf m}$, что {при всех ${\bf m}\in{\mathbb Z}_+^\kappa\setminus\{0\}$} выполняется неравенство
\begin{equation}\label{th_eq3}
\|\Delta^{\ell}C(t)\|_R\leqslant \alpha\left\|L(\lambda_m+\langle{\bf m},{\bf r}\rangle+d/dt)C(t)\right\|_R.
\end{equation}

Все необходимые вспомогательные оценки сделаны. Вернёмся к отображению $h'_\eta(0,0)$ из пространства $H^{\ell}$ в пространство $H^0$, то есть
$$
h'_\eta(0,0):\eta_1:=\sum_{{\bf m}\in\mathbb{Z}_+^{\kappa}\setminus\{0\}}a_{\bf m}(t){\bf x}^{\bf m}\mapsto\sum_{{\bf m}\in\mathbb{Z}_+^{\kappa}\setminus\{0\}}
L\Bigl(\lambda_m+\langle{\bf m},{\bf r}\rangle+\frac{d}{dt}\Bigr)a_{\bf m}(t){\bf x}^{\bf m}=:\eta_2.
$$
Инъективность и сюръективность $h'_\eta(0,0)$ следуют из того, что линейное дифференциальное уравнение
$$
L\Bigl(\lambda_m+\langle{\bf m},{\bf r}\rangle+\frac{d}{dt}\Bigr)y=0
$$
не имеет ненулевых полиномиальных решений, а уравнение
$$
L\Bigl(\lambda_m+\langle{\bf m},{\bf r}\rangle+\frac{d}{dt}\Bigr)y=P(t)
$$
обладает полиномиальным решением той же степени, что и полином $P$.

Наконец, проверим непрерывность обратного отображения $h'_\eta(0,0)^{-1}$ из пространства $H^0$ в пространство $H^{\ell}$,
$$
h'_\eta(0,0)^{-1}:\eta_3:=\sum_{{\bf m}\in\mathbb{Z}_+^{\kappa}\setminus\{0\}}b_{\bf m}(t)\,{\bf x}^{\bf m}\mapsto
\sum_{{\bf m}\in\mathbb{Z}_+^{\kappa}\setminus\{0\}}\tilde b_{\bf m}(t)\,{\bf x}^{\bf m}=:\eta_4,
$$
где $L(\lambda_m+\langle{\bf m},{\bf r}\rangle+d/dt)\tilde b_{\bf m}=b_{\bf m}$.
По определению норм пространств $H^0$ и $H^{\ell}$, с учётом неравенств \eqref{th_eq2} при $R>\mathcal{K}\theta$ и \eqref{th_eq3} получаем, что
$$
\|\eta_4\|_\ell=\sum_{{\bf m}\in\mathbb{Z}_+^{\kappa}\setminus\{0\}}
\frac{\displaystyle(|\lambda_m+\langle{\bf m},{\bf r}\rangle|+\mathcal{K}|{\bf m}|)^{\ell}\,\|\tilde{b}_{\bf m}(t)\|_R} {\left|\Gamma(\langle{\bf m},{\bf r}\rangle/s)\right|}\leqslant \Upsilon\sum_{{\bf m}\in\mathbb{Z}_+^{\kappa}\setminus\{0\}}
\frac{\displaystyle\,\|\Delta^{\ell}\,\tilde{b}_{\bf m}(t)\|_R}{\left|\Gamma(\langle{\bf m},{\bf r}\rangle/s)\right|}\leqslant
$$
$$
\leqslant\Upsilon\alpha\sum_{{\bf m}\in\mathbb{Z}_+^{\kappa}\setminus\{0\}}
\frac{\displaystyle\,\|b_{\bf m}(t)\|_R}{\left|\Gamma(\langle{\bf m},{\bf r}\rangle/s)\right|}=\Upsilon\alpha\|\eta_3\|_0,
$$
где $\displaystyle\Upsilon=\Theta^{\ell}\sup_{\bf m\in\mathbb{Z}_+^{\kappa}\setminus \{0\}}\left(\left(|\lambda_m+\langle{\bf m},{\bf r}\rangle|+\mathcal{K}|
{\bf m}|\right)^{\ell}/|{\bf m}|^{\ell}\right)$.
\end{itemize}
\medskip

Таким образом, применение к $h(\rho,\eta)$ теоремы о неявном отображении правомерно и, следовательно, согласно этой теореме существуют $\rho>0$ и $\eta_\rho\in H^\ell$ такие, что
$$
L(\lambda_m+\hat{\delta})\eta_{\rho}+\widehat{L}(t,\rho \,{\bf x},\lambda_m+\hat{\delta})\eta_{\rho}+
\widehat N(t,\rho\, {\bf x},\,\rho^{|{\boldsymbol \tau}|}{\bf x}^{\boldsymbol \tau}\,\eta_{\rho},\ldots,
\rho^{|{\boldsymbol \tau}|}{\bf x}^{\boldsymbol \tau}\,(\lambda_m+\hat{\delta})^n\eta_{\rho})=0.
$$
Поскольку замена переменных ${\bf x}\mapsto{\bf x}/\rho$ порождает автоморфизм
${\mathbb C}[t][[x_1,\ldots,x_\kappa]]_*\ni\eta({\bf x})\mapsto\eta({\bf x}/\rho)\in{\mathbb C}[t][[x_1,\ldots,x_\kappa]]_*$, коммутирующий с оператором $\hat{\delta}$, мы можем заметить, что степенной ряд $\eta_{\rho}\left({\bf x}/\rho\right)$, так же как и ряд \eqref{m_resh}, удовлетворяет уравнению \eqref{m_eq}.
Следовательно, эти два ряда совпадают между собой в силу единственности (см. замечание 1) и тогда {$\hat\psi\in H^\ell$ и} ряд
$$
\sum_{{\bf m}\in\mathbb{Z}_+^{\kappa}\setminus\{0\}}\frac{\|c_{\bf m}\|_R}{|\Gamma(\langle{\bf m},{\bf r}\rangle/s)|}{\bf x}^{\bf m}
$$
имеет ненулевой радиус сходимости. { Допустим, он сходится при $|x_i|\leqslant\varrho$, $i=1,\ldots,\kappa$. Остаётся показать равномерную сходимость ряда $$
\sum_{{\bf m}\in\mathbb{Z}_+^{\kappa}\setminus\{0\}}\frac{c_{\bf m}(t)}{\Gamma(\langle{\bf m},{\bf r}\rangle/s)}x^{\langle{\bf m},{\bf r}\rangle}=\sum_{k=m+1}^{\infty}\frac{c_k(t)}{\Gamma((\lambda_k-\lambda_{m})/s)}x^{\lambda_k-\lambda_{m}}
$$
во всяком открытом секторе $S$ с вершиной в нуле достаточно малого радиуса и раствора меньше $2\pi$.

Воспользуемся оценками $|x^{r_i}|\leqslant \gamma\,|x|^{{\rm Re}\,r_i}$ и ${\rm Re}\langle{\bf m},{\bf r}\rangle\geqslant\beta\,|{\bf m}|$ при подходящих фиксированных $\beta,\gamma>0$. Зафиксируем $\epsilon<\beta/\mathcal{K}$ и пусть радиус сектора $S$ настолько мал, что выполнено
$$
|t|=|\ln x|\leqslant |x|^{-\epsilon} \qquad \forall x\in S.
$$
Тогда
$$
|c_{\bf m}(t)|\leqslant \|c_{\bf m}\|_R\,|t|^{\deg c_{\bf m}}\leqslant \|c_{\bf m}\|_R \,|x|^{-\epsilon\,\deg c_{\bf m}}\leqslant \|c_{\bf m}\|_R\,
|x|^{-\epsilon\,\mathcal{K}|{\bf m}|}
$$
и, следовательно,
$$
|c_{\bf m}(t)x^{\langle{\bf m},{\bf r}\rangle}|=|c_{\bf m}(t)|\,|x^{r_1}|^{m_1}\ldots|x^{r_{\kappa}}|^{m_\kappa}\leqslant
\|c_{\bf m}\|_R\,\gamma^{|{\bf m}|}|x|^{{\rm Re}\langle{\bf m},{\bf r}\rangle}\,|x|^{-\epsilon\,\mathcal{K}|{\bf m}|}\leqslant
$$
$$
\leqslant\|c_{\bf m}\|_R\,\gamma^{|{\bf m}|}|x|^{(\beta-\epsilon\,\mathcal{K})|{\bf m}|}\leqslant\|c_{\bf m}\|_R\,\varrho^{|{\bf m}|},
$$
если $\gamma\,|x|^{\beta-\epsilon\,\mathcal{K}}\leqslant\varrho$, откуда следует нужная сходимость.}


\begin{thebibliography}{99}

\bibitem{Du} H.~Dulac, Sur les cycles limites, {\it Bull. Soc. Math. France}, {\bf 51} (1923), 45--188.

\bibitem{GG_D}
Р.\,Р.~Гонцов, И.\,В.~Горючкина, О сходимости формальных рядов Дюлака, удовлетворяющих алгебраическому ОДУ, {\it Матем. сб.}, {\bf 210}:9 (2019), 3--18.

\bibitem{b}
А.\,Д.~Брюно, Асимптотики и разложения решений обыкновенного дифференциального уравнения, {\it Успехи матем. наук}, {\bf 59}:3 (2004), 31--80.

\bibitem{RamisSibuya} J.-P. Ramis, Y. Sibuya, Hukuhara's domains and fundamental existence and uniqueness theorems for asymptotic solutions of Gevrey type,
{\it Asympt. Anal.}, {\bf 2}:1 (1989), 39--94.

\bibitem{Mal}
B.\,Malgrange, Sur le th\'eor\`eme de Maillet, {\it Asympt. Anal.}, {\bf 2}:1 (1989), 1--4.

\bibitem{GG_MM} R. Gontsov, I. Goryuchkina, The Maillet--Malgrange type theorem for generalized power series, {\it Manuscr. Math.}, {\bf 156} (2018), 171--185.

\bibitem{GGarxiv}
Р.\,Р.\,Гонцов, И.\,В.\,Горючкина,  Сходимость обобщённых степенных рядов, удовлетворяющих функциональным уравнениям, {\it Успехи матем. наук},
{\bf 80}:3 (2025), 3--66.

\bibitem{Dick}
L.\,E.\,Dickson, Finiteness of the odd perfect and primitive abundant numbers with $n$ distinct prime factors, {\it Amer. J. Math.}, \textbf{35}:4 (1913), 413--422.

\bibitem{Batm} Г.\,Бейтман, А.\,Эрдейи, Высшие трансцендентные функции, Т. 1,  М.: Наука, 1973.

\bibitem{KolmFomin}
А.\,Н.\,Колмогоров, С.\,В.\,Фомин, Элементы теории функций и функционального анализа, 7-е изд., М.: Физматлит, 2006.

\end{thebibliography}
\end{document}